# A Congruence Property of Solvable Polynomials

Nicholas Phat Nguyen[1]

**Abstract.**  We describe a congruence property of solvable polynomials over **Q**, based on the irreducibility of cyclotomic polynomials over number fields that meet certain conditions.[2]

**A.      INTRODUCTION.**   Consider a monic polynomial $h(X) \in \mathbf{Z}[X]$ with integer coefficients. A landmark result in the history of mathematics says that $h(X)$ is solvable by radicals if and only if the splitting extension generated by $h(X)$ has a solvable Galois group over **Q**. In this note, we will describe a congruence property of solvable polynomials.

Specifically, we will prove the following.

<u>**Theorem**</u> :  *Let $h(X) \in \mathbf{Z}[X]$ be a monic polynomial with integer coefficients, and let P be the set of rational prime numbers p such that $h(X)$ splits completely into distinct linear factors mod p. If the Galois group of $h(X)$ has a nontrivial solvable quotient group, then there is an integer n having a nontrivial common factor with the discriminant of $h(X)$, such that every generating set of the unit group $(\mathbf{Z}/n\mathbf{Z})^x$ has at least a residue class containing <u>no</u> rational prime number p in P.*

The theorem applies to polynomials solvable by radicals because their Galois groups are solvable. But it also applies to irreducible polynomials $h(X)$ whose Galois group is the

---

full symmetric group on deg($h$) letters, such as the case with many polynomials of the type $X^m + aX + b$, because those Galois groups have normal subgroups of index 2.

This theorem is not constructive because the number $n$ is unspecified, although we know $n$ must have a nontrivial common factor with the discriminant of $h(X)$. Nonetheless, the theorem gives us a property of solvable polynomials that depends only on the discriminant of the polynomial and its factorization modulo rational primes, rather than on the roots of the polynomial and their Galois property.

Let $K$ be the splitting extension of $h(X)$. The set $P$ has a density of $1/\deg K$ according to the Chebotarev density theorem. The theorem says that the distribution of $P$ must have a certain gap modulo some number $n$. What we mean is that every generating set of the unit group $(\mathbb{Z}/n\mathbb{Z})^\times$ must include at least a residue class containing no prime numbers in $P$.

For polynomials $h(X)$ where the corresponding set $P$ can be characterized by some congruence conditions, such a distribution gap modulo some number $n$ is easy to see.

- For example, if $h(X) = X^q - r$ with $q > 2$, then $h(X)$ is obviously solvable and its splitting extension contains $\mathbb{Q}(\zeta_q)$. For $h(X)$ to split completely into distinct linear factors modulo a rational prime $p$, $p$ must be $\equiv 1 \bmod q$, as the multiplicative group $(\mathbb{Z}/p\mathbb{Z})^\times$ must contain elements whose order is a multiple of $q$. Any generating set of $(\mathbb{Z}/q\mathbb{Z})^\times$ must of course include a residue class that is $\not\equiv 1 \bmod q$, and such a residue class has no prime $p$ such that $h(X)$ splits completely modulo $p$.
- Consider an irreducible quadratic polynomial $h(X)$. To the quadratic extension $E$ defined by $h(X)$ we can associate a real primitive character $\chi$: $(\mathbb{Z}/c\mathbb{Z})^\times \to \{\pm 1\}$, where $c$ is equal to the discriminant of $E$, such that for any rational prime $p$ not dividing $c$, $\chi(p) = 1$ if and only if $p$ splits completely in $E$. Those rational prime numbers include the set $P$, so that $\chi(p) = 1$ for any prime number $p$ in the set $P$. Any generating set of $(\mathbb{Z}/c\mathbb{Z})^\times$ must of course contain residue classes $a \bmod c$ such that $\chi(a) = -1$ because the character $\chi$ is primitive. Such a residue class contains no rational prime number in $P$.

However, if we do not have a congruence condition to describe the set $P$ associated with a given polynomial $h(X)$, such as the case for irreducible cubic and quartic polynomials whose Galois group is non-abelian or the case for polynomials of the type $X^m + aX + b$ whose Galois group is the full symmetric group on $m$ letters, the congruence property described in the theorem could perhaps provide a helpful perspective.

B. **SOLVABLE EXTENSIONS OF Q**. We begin with some notation. Let $K$ be a finite extension of **Q** (a number field) with ring of integers $A$. The ring $A$ is a Dedekind domain and any integral ideal $pA$ generated by a rational prime $p$ has a unique prime ideal factorization $\prod m_j$ (with possible repeated factors), where each $m_j$ is a maximal ideal in $A$ with $m_j \cap \mathbf{Z} = p\mathbf{Z}$. The residue field $A/m_j$ is an extension of the finite field $\mathbf{F}_p = \mathbf{Z}/p\mathbf{Z}$, and the degree of that residue extension is known as the residue degree of $m_j$ over the prime $p$ or over **Q**.

Of particular interest to us are the rational primes $p$ such that at least one of the maximal ideals $m_j$ in the prime ideal factorization of $pA$ has residue degree one over $p$, that is to say the residue field $A/m_j$ is isomorphic to $\mathbf{F}_p$. These rational primes do not have a particular name in the literature, but for convenient reference we will refer to them in this note as primes that are *semi-split* in $K$. (Rational primes $p$ which are unramified in $K$ and such that all maximal ideals in the factorization of $pA$ have residue degree one are said to split totally or completely in $K$.)

If $K = \mathbf{Q}(\alpha)$ with an integral element $\alpha$ whose minimal polynomial is $f(X) \in \mathbf{Z}[X]$, then with a finite number of exceptions, the rational primes $p$ semi-split in K are precisely the rational primes $p$ such that $f(X) \bmod p$ has a root in $\mathbf{Z}/p\mathbf{Z}$, i.e., the prime factors $p$ of the numbers $f(k)$ as $k$ runs through **Z**.

Our theorem depends on the following necessary condition for a number field $K$ to contain a nontrivial solvable extension of **Q**.

<u>Proposition 1</u>.  *Let K be a number field. For K to contain a nontrivial solvable extension of **Q**, it is necessary that we can find an integer n where every generating set of $(\mathbf{Z}/n\mathbf{Z})^x$ has at least a residue class containing <u>no</u> rational prime number semi-split in K.*

*Proof.* By the definition of solvable extensions, $K$ contains a nontrivial solvable extension of $\mathbf{Q}$ if and only if $K$ contains a nontrivial abelian extension of $\mathbf{Q}$. According to the Kronecker-Weber theorem, any abelian extension of $\mathbf{Q}$ is contained in some cyclotomic extension $\mathbf{Q}(\zeta_n)$. So $K$ contains a nontrivial abelian extension of $\mathbf{Q}$ if and only if $K \cap \mathbf{Q}(\zeta_n) \neq \mathbf{Q}$ for some $n$.

For any integer $n > 0$, we define the $n^{\text{th}}$-cyclotomic polynomial $\Phi_n(X)$ as $\prod (X - \zeta)$, with the product running over all primitive $n^{\text{th}}$-roots of unity $\zeta$. The coefficients of $\Phi_n(X)$ are rational numbers that are integral over $\mathbf{Z}$, and so is $\Phi_n(X)$ a monic polynomial in $\mathbf{Z}[X]$. It is well-known that $\Phi_n(X)$ is irreducible over $\mathbf{Q}$, so that $[\mathbf{Q}(\zeta_n) : \mathbf{Q}] =$ degree of $\Phi_n(X)$.

To say that $K \cap \mathbf{Q}(\zeta_n) \neq \mathbf{Q}$ is equivalent to saying that $\Phi_n(X)$ is reducible over $K$. Indeed, for any extension $K$ of $\mathbf{Q}$, $\Phi_n(X)$ is irreducible over $K$ if and only if we have $[K(\zeta_n) : K] = [\mathbf{Q}(\zeta_n) : \mathbf{Q}]$, which is the case if and only if $\mathbf{Q}(\zeta_n) \cap K = \mathbf{Q}$ because $\mathbf{Q}(\zeta_n)$ is a Galois extension of $\mathbf{Q}$.

For $\Phi_n(X)$ to be reducible over $K$, it is necessary that the number field $K$ does not meet the condition of Proposition 2 below, which gives us a sufficient condition for $\Phi_n(X)$ to be irreducible over $K$. In particular, if a number field $K$ fails the condition of Proposition 2 for $n$, every generating set of the group $(\mathbf{Z}/n\mathbf{Z})^\times$ must include at least a residue class containing no rational prime number semi-split in $K$. ∎

**C. CYCLOTOMIC POLYNOMIALS REVISITED.** Gauss was the first person to prove that $\Phi_n(X)$ is irreducible over $\mathbf{Q}$ when $n$ is prime. Many other distinguished mathematicians have also come up with proofs for the irreducibility of $\Phi_n(X)$ over $\mathbf{Q}$, including Schönemann and Eisenstein (for $n$ prime), Kronecker (for $n$ prime and for general $n$), Dedekind, Landau, and Schur (for general $n$). For an outline of some of these proofs, the reader can see [3] and the references cited in that article. Dedekind himself gave three proofs for the irreducibility of $\Phi_n(X)$ over $\mathbf{Q}$, which tells us that Dedekind regarded the irreducibility of cyclotomic polynomials as important. In fact, there is a connection that can be drawn between the irreducibility of cyclotomic polynomials and the Artin reciprocity homomorphism at the foundation of class field theory.

We will generalize a proof by Dedekind to prove the irreducibility of $\Phi_n(X)$ over number fields that meet certain conditions.

Incidentally, Kronecker was the first person who thought of proving the irreducibility of $\Phi_n(X)$ over number fields other than $\mathbf{Q}$. See [4]. Kronecker discovered that he could prove irreducibility of the cyclotomic polynomial in the case of general $n$ if he could show that when $n$ is a prime power, $\Phi_n(X)$ is irreducible over the cyclotomic field $\mathbf{Q}(\zeta_m)$ if $m$ is relatively prime to $n$. In modern language, his approach relies on the essential fact that if $n$ is a power of a prime $p$, then $p$ ramifies completely in the extension $\mathbf{Q}(\zeta_n)$, while it is unramified in $\mathbf{Q}(\zeta_m)$.

More generally, a common sufficient condition for $\Phi_n(X)$ to be irreducible over a number field $K$ is that $n$ is relatively prime to the discriminant of $K$. We provide below another sufficient condition.

**Proposition 2:** *Let $K$ be a number field. If each residue class in a generating set of the group $(\mathbf{Z}/n\mathbf{Z})^{\times}$ contains a rational prime semi-split in $K$, then the $n^{th}$ cyclotomic polynomial $\Phi_n(X)$ is irreducible over $K$.*

*Proof.* Note that each irreducible polynomial in $K[X]$ corresponds to an orbit in a given separable closure of $K$ under the action of the absolute Galois group of $K$. All elements in that orbit are the roots of the irreducible polynomial.

If $\Phi_n(X)$ is not irreducible over $K$, then there must be two different roots of $\Phi_n(X)$, say $u$ and $v$, that belong to two distinct Galois orbits over $K$. The integral elements $u$ and $v$ must then have distinct minimal polynomials $g$ and $h$ over $K$ whose product divides $\Phi_n(X)$ in $A[X]$.

Let $B$ be the ring of integers in the extension $K(u, v)$ of $K$. For each maximal ideal $\mathfrak{m}$ of $A$ such that $\Phi_n(X)$ mod $\mathfrak{m}$ is separable, consider any maximal ideal $\wp$ of $B$ in the prime ideal factorization of $\mathfrak{m}B$. The images $u$ mod $\wp$ and $v$ mod $\wp$ in the residue field $B/\wp$ must belong to different Galois orbits over $A/\mathfrak{m}$, because they are the roots of polynomials $g$ mod $\mathfrak{m}$ and $h$ mod $\mathfrak{m}$ whose product divides the separable polynomial $\Phi_n(X)$ mod $\mathfrak{m}$.

Because both $u$ and $v$ are primitive $n^{th}$ roots of unity, there is an invertible residue class $a$ mod $n$ such that $u^a = v$. For our purpose, we can assume that $a$ mod $n$ belongs to the given generating set for the group $(\mathbf{Z}/n\mathbf{Z})^\times$. Indeed, if both $u$ and $u^a$ belong to the same Galois orbit over $K$ regardless of which primitive $n^{th}$ root of unity $u$ we choose and which residue class $a$ mod $n$ in the generating set we choose, then all the primitive $n^{th}$ roots of unity must be in the same Galois orbit over $K$, contrary to our initial assumption.

By hypothesis, we can find a rational prime $p$ that is semi-split in $K$ and congruent to $a$ mod $n$. This rational prime $p$ does not divide $n$, and hence the polynomial $X^n - 1$ is separable over any field of characteristic $p$. In particular, for each maximal ideal $\mathfrak{m}$ in $A$ such that $\mathfrak{m} \cap \mathbf{Z} = p\mathbf{Z}$, the polynomial $X^n - 1$ is separable over the residue field $A/\mathfrak{m}$. That means $\Phi_n(X)$ is also separable over $A/\mathfrak{m}$. In light of our foregoing discussion, if $\wp$ is any maximal ideal of $B$ in the prime ideal factorization of $\mathfrak{m}B$, then the elements $u$ mod $\wp$ and $v$ mod $\wp$ will be in distinct Galois orbits over $A/\mathfrak{m}$.

The congruence condition on $p$ means that $v$ is equal to the $p^{th}$ power of $u$, so that the element $v$ mod $\wp$ is the transform of $u$ mod $\wp$ under the Frobenius automorphism which raises each element in a field of characteristic $p$ to the $p^{th}$ power. Note that each element algebraic over $\mathbf{F}_p$ and its $p^{th}$ power must belong to the same Galois orbit over $\mathbf{F}_p$, because the Frobenus automorphism is a generator of any finite Galois group over $\mathbf{F}_p$.

By the semi-split condition for $p$, we can choose a maximal ideal $\mathfrak{m}$ of $A$ sitting over $p$ such that $A/\mathfrak{m}$ is isomorphic to $\mathbf{F}_p$. That means $u$ mod $\wp$ and $v$ mod $\wp$ must belong to the same Galois orbit over $A/\mathfrak{m}$. However, we saw earlier that $u$ mod $\wp$ and $v$ mod $\wp$ must be in different Galois orbits over $A/\mathfrak{m}$, in part because $u$ and $v$ are assumed to be in different Galois orbits over $K$. This contradiction shows that all the roots of $\Phi_n(X)$ must be in one Galois orbit over $K$. Therefore $\Phi_n(X)$ must be irreducible over $K$. ∎

Let us illustrate Proposition 2 by a classical example. Consider the cyclotomic extension $\mathbf{Q}(\zeta_m)$ for any integer $m$ relatively prime to $n$. The rational primes that are unramified and semi-split in $\mathbf{Q}(\zeta_m)$ are the primes $p \equiv 1 \pmod{m}$. The Chinese remainder theorem and Dirichlet's theorem combine to tell us that there are always rational primes $p$ (indeed infinitely many) that satisfy both the congruence $p \equiv 1 \pmod{m}$ and the congruence $p \equiv a \pmod{n}$ for each invertible residue class $a$ mod $n$. Each such rational prime would be semi-split in $\mathbf{Q}(\zeta_m)$ and would also belong to the residue class ($a$ mod $n$). Accordingly, the $n^{\text{th}}$ cyclotomic polynomial $\Phi_n(X)$ is irreducible over $\mathbf{Q}(\zeta_m)$.

Note that Proposition 2 implies immediately that if all but a finite number of rational primes in $\mathbf{Q}$ are semi-split in a number field $K$, then all cyclotomic polynomials $\Phi_n(X)$ are irreducible over $K$, and $K$ cannot contain a nontrivial solvable extension of $\mathbf{Q}$. That is because by Dirichlet's theorem, each invertible residue class $a$ mod $n$ (for any $n$) contains infinitely many prime numbers. So if all but a finite number of rational primes in $\mathbf{Q}$ are semi-split in $K$, then each invertible residue class $a$ mod $n$ must contain infinitely many primes semi-split in $K$, making $\Phi_n(X)$ irreducible over $K$ for any $n$. In fact, $K = \mathbf{Q}$ in this case, because the only number field in which the semi-split rational primes have density one is $\mathbf{Q}$ itself. See [2], Exercise 6.2.

D. **PROOF OF THEOREM.** Let us now apply the result of Proposition 1 to a monic polynomial $h(X) \in \mathbf{Z}[X]$ with integer coefficients, with splitting extension $K$. If the Galois group of $h(X)$ has a nontrivial solvable quotient group, that means $K$ contains a nontrivial solvable extension of $\mathbf{Q}$. According to Proposition 1, there must be an integer $n$ such that every generating set of $(\mathbf{Z}/n\mathbf{Z})^\times$ has at least a residue class containing <u>no</u> rational prime number semi-split in $K$. In particular, such a residue class contains no rational prime number that is unramified and semi-split in $K$.

Because $K$ is a Galois extension, a rational prime number is unramified and semi-split in $K$ if and only if it splits completely in $K$. Moreover, because $K$ is generated by the splitting fields for each irreducible factor of the polynomial $h(X)$, a rational prime number splits completely in $K$ if and only if it splits completely in each of those splitting fields.

If $g(X)$ is an irreducible factor of $h(X)$, and $L = \mathbf{Q}(\alpha)$ is the simple extension obtained by adjoining a root $\alpha$ of $g(X)$ to $\mathbf{Q}$, then a rational prime number splits completely in the splitting field of $g(X)$ if and only if it splits completely in $L$. See, e.g., [1] at Section 16.

The rational primes $p$ such that $g(X)$ is split completely into distinct linear factors modulo $p$ are precisely the rational prime numbers $p$ that are: (a) relatively prime to the algebraic integer $g'(\alpha)$, and (b) split completely in $L$. Indeed, if a rational prime $p$ is relatively prime to the algebraic integer $g'(\alpha)$, then the prime ideal decomposition of $p$ in $L$ is determined by the factorization of the polynomial $g(X)$ modulo $p$. See, e.g., Theorem 14.10 of [1]. In that case, the rational prime number $p$ splits completely in $L$ if and only if $g(X)$ splits completely into linear factors modulo $p$. These linear factors modulo $p$ are necessarily distinct because $p$ is relatively prime to $g'(\alpha)$. Conversely, if $g(X)$ splits completely into distinct linear factors modulo $p$, then $p$ must be relatively prime to $g'(\alpha)$. The factorization of $g(X)$ into distinct linear factors modulo $p$ then translates into the complete factorization of the ideal generated by $p$ into distinct maximal ideals all having residue degree 1.

To summarize, if the splitting extension $K$ contains a nontrivial solvable extension of $\mathbb{Q}$, then there is an integer $n$ such that every generating set of the unit group $(\mathbb{Z}/n\mathbb{Z})^{\times}$ has at least a residue class containing no rational prime number semi-split in $K$, which implies that such a residue class contains no rational prime $p$ such that $h(X)$ splits completely into distinct linear factors modulo $p$.

While our proof does not determine the number $n$, it tells us that $K \cap \mathbb{Q}(\zeta_n) \neq \mathbb{Q}$. The rational primes that ramify in $K \cap \mathbb{Q}(\zeta_n)$ must divide $n$ because only those rational prime numbers ramify in $\mathbb{Q}(\zeta_n)$. At the same time, any rational prime that ramifies in the extension $K \cap \mathbb{Q}(\zeta_n)$ must ramify in at least one simple extension obtained by adjoining a root of $h(X)$ to $\mathbb{Q}$. Otherwise such a rational prime would be unramified in all the splitting fields for each irreducible factor of the polynomial $h(X)$, and therefore also unramified in $K$.

If we again let $g(X)$ be an irreducible factor of $h(X)$, and $L = \mathbb{Q}(\alpha)$ the simple extension obtained by adjoining a root $\alpha$ of $g(X)$ to $\mathbb{Q}$, then a rational prime number $p$ ramifies in $L$ only if it divides the algebraic integer $g'(\alpha)$, because otherwise the prime ideal factorization of $p$ in $L$ would mirror the factorization of $g(X)$ modulo $p$, and $g(X)$ would be separable modulo $p$ if $p$ does not divide $g'(\alpha)$. If $p$ divides $g'(\alpha)$, then $p$ also divides $h'(\alpha)$ and therefore the discriminant of $h(X)$ as well.

Accordingly, any rational prime $p$ that ramifies in $K \cap \mathbf{Q}(\zeta_n)$ must divide $n$ as well as the discriminant of $h(X)$. Because $K \cap \mathbf{Q}(\zeta_n)$ is a nontrivial extension of $\mathbf{Q}$, there are rational prime numbers that ramify in $K \cap \mathbf{Q}(\zeta_n)$, and therefore the number $n$ has a nontrivial common factor with the discriminant of $h(X)$.

The proof of our theorem is now complete. ∎

**Acknowledgments by Author**: I am very grateful to Keith Conrad for his helpful comments on earlier drafts of this paper.

## REFERENCES

[1] Goro Shimura (2010). *Arithmetic of Quadratic Forms.* Springer Monograph in Mathematics.

[2] J.W.S. Cassels & A. Frohlich (Second Edition 2010). *Algebraic Number Theory.* London Mathematical Society.

[3] Steven H. Weintraub (2013). *Several Proofs of the Irreducibility of the Cyclotomic Polynomials*. The American Mathematical Monthly, 120:6, 537-545.

[4] Leopold Kronecker (1854). *Memoire sur les facteurs irreductibles de l'expression $x^n - 1$*. Journal de mathematiques pures et appliquees, 1re serie, tome 19, 177-192.